\documentclass{article}

\begin{document}

\title{Adiabatic limits of co-associative Kovalev-Lefschetz fibrations}
\author{Simon Donaldson}
\maketitle


\newtheorem{lem}{Lemma}
\newtheorem{prop}{Proposition}
\newtheorem{cor}{Corollary}
\newtheorem{defn}{Definition}
\newtheorem{conj}{Conjecture}
\newcommand{\bP}{{\bf P}}
\newcommand{\bC}{{\bf C}}
\newcommand{\hook}{\rightharpoonup}
\newcommand{\uomega}{\underline{\omega}}
\newcommand{\uTheta}{\underline{\Theta}}
\newcommand{\umu}{\underline{\mu}}
\newcommand{\ulambda}{\underline{\lambda}}
\newcommand{\tuomega}{\tilde{\underline{\omega}}}
\newcommand{\tuTheta}{\tilde{\underline{\Theta}}}
\newcommand{\tumu}{\tilde{\underline{\mu}}}
\newcommand{\tulambda}{\tilde{\underline{\lambda}}}
\newcommand{\bR}{{\bf R}}
\newcommand{\bZ}{{\bf Z}}
\newcommand{\ut}{{\underline{t}}
\newcommand{\Vol}{\rm Vol}}

\ \ \ \ \ \ \ \ \ \ \ \ \ \ \ \ \ \ \ {\it To Maxim Kontsevich, on his 50th. birthday}

\

\

\

\section{Introduction} This article makes the first steps in what we hope will be a longer project, investigating seven-dimensional $G_{2}$-manifolds from the point of view of co-associative fibrations, and in particular the \lq\lq adiabatic limit'', when the diameters of the fibres shrink to zero.   To set up the background, recall that there is a notion of a \lq\lq positive'' exterior 3-form on an oriented $7$-dimensional real vector space $E$, and that such a form defines a Euclidean structure on $E$ (definitions are given in Section 2 below). Thus  there is a notion of a positive $3$-form  $\phi$ on an oriented $7$-manifold $M$, which defines a Riemannian metric $g_{\phi}$ and hence a $4$-form $*_{\phi} \phi$ (writing $*_{\phi}$ for the Hodge $*$-operator of $g_{\phi}$).  A {\it torsion-free $G_{2}$-structure} can be defined to be a positive form $\phi$ such that $\phi$ and $*_{\phi}\phi$ are both closed. (We refer to \cite{kn:J} for the terminology and background.)  We can also consider the weaker condition of a {\it closed $G_{2}$-structure} in which we just require that the $3$-form $\phi$ is closed. We have then various fundamental questions, such as the following.
\begin{enumerate}
 \item Given a compact oriented $7$-manifold and a de Rham cohomology class $C\subset \Omega^{3}(M)$, does $C$ contain a positive form? In other words, is there a non-empty  subset $C^{+}\subset C$ of positive forms?
\item   If so, is there a torsion free $G_{2}$-structure in $C^{+}$, i.e. can we solve the equation $d *_{\phi} \phi =0$ for $\phi\in C^{+}$?
\item The {\it Hitchin functional} $H$ on the space of positive 3-forms is just the volume
$$    H(\phi) = {\rm Vol} (M, g_{\phi}). $$
The equation $d *_{\phi} \phi =0$ is the Euler-Lagrange equation for this functional restricted to $C^{+}$---in fact critical points are local maxima \cite{kn:H1}. It is then interesting to ask if $H$ is bounded above on $C^{+}$ and, if a critical  point exists, whether it is a global maximum. 
\item Bryant \cite{kn:Br} defined a \lq\lq Laplacian'' flow on $C^{+}$:
$$   \frac {\partial \phi}{\partial t} = d \left(*_{\phi} d *_{\phi} \phi\right), $$
whose fixed points are torsion-free $G_{2}$ structures. It is interesting to ask if this flow converges as $t\rightarrow \infty$ to a fixed point.   \end{enumerate}
 While  these are obvious and natural questions, any kind of systematic answers seem far out of reach at present. In this article we will study 7-manifolds with extra structure given by a \lq\lq Kovalev-Lefschetz fibration'' $\pi:M\rightarrow S^{3}$ with co-associative fibres. These fibres will be $K3$ surfaces, either smooth or with nodal singularities. The idea underlying  our discussion is that there should be an \lq\lq adiabatic limit'' of the torsion-free $G_{2}$-condition, and more generally of each of the questions above. Leaving for the future any precise statement along those lines we will write down precise definitions which we propose as limiting  geometric objects and we will see that questions (1)-(4) have sensible analogs in this setting.

Co-associative fibrations were studied by Gukov, Yau and Zaslow \cite{kn:GYZ} and there is overlap between our discussion in Section 2 and theirs. More recent work by Baraglia \cite{kn:B} in the case of torus fibrations is  especially relevant to our constructions in this paper (see also the brief comments in 5.1 below) In particular Baraglia made the connection with   maximal submanifolds in spaces with indefinite signature, which is the fundamental idea that we use here.  The notion of a  Kovalev-Lefschetz fibration was essentially introduced by Kovalev in \cite{kn:K}. In other directions there is a substantial literature on adiabatic limits of various kinds of structures in Riemannian geometry.  Examples include  constant scalar curvature K\"ahler metrics  \cite{kn:F} and many papers by different authors on \lq\lq collapsing'' Calabi-Yau manifolds, such as \cite{kn:GW}, \cite{kn:T}.
We will not attempt to give a survey of existence results for compact $G_{2}$-manifolds here,  but we mention the landmark results of Joyce \cite{kn:J0}, Kovalev \cite{kn:K0} and Corti, Haskins, Nordstr\"om, Pacini \cite{kn:CHNP}.  
\section{Co-associative fibrations}

We review some standard multi-linear algebra background. Let $E$ be an oriented $7$-dimensional real vector space and let $\phi$ be a 3-form in $\Lambda^{3} E^{*}$. We have then a quadratic form on $E$ with values in the oriented line $\Lambda^{7} E^{*}$  defined by 
\begin{equation} G_{\phi}(v) = (v\hook \phi)^{2} \wedge \phi. \end{equation}
We say that $\phi$ is a {\em positive} $3$-form if $G_{\phi}$ is positive definite. In this case, any choice of (oriented) volume form makes $G_{\phi}$ into a Euclidean structure on $E$ and there is a unique choice of the volume form such that $\vert \phi\vert^{2} =7$. So we now have a Euclidean structure $G_{\phi}$ canonically determined by a positive 3-form $\phi$. 

For $(E,\phi)$ as above we say that a 4-dimensional subspace $V\subset E$ is {\it co-associative } if the restriction of $\phi$ to $V$ vanishes. A standard model for such a triple $(E,\phi, V)$ can be built as follows. We start with $\bR^{4}$ with co-ordinates $x_{0}, \dots, x_{3}$ and standard metric and orientation. Let $\Lambda^{2}_{+}$ be the space of {\it self-dual 2-forms}. This has a standard basis $\omega_{1}, \omega_{2}, \omega_{3}$ given by  $$  \omega_{i}= dx_{0}dx_{i} + dx_{j} dx_{k} $$
where $i,j,k$ run over cyclic permutations of $1,2,3$. Let $t_{1}, t_{2}, t_{3}$ be the co-ordinates on $\Lambda^{2}_{+}$ corresponding to this basis and set
\begin{equation} \phi_{0}= \sum_{i} \omega_{i} dt_{i} - dt_{1} dt_{2} dt_{3}, \end{equation}
an element of $\Lambda^{3} (E_{0})^{*}$ where $E_{0}= \bR^{4}\oplus (\Lambda^{2}_{+})^{*}$. This is a positive $3$-form if $E_{0}$ is given the orientation
\begin{equation} - dt_{1} dt_{2} dt_{3} dx_{0} dx_{1} dx_{2} dx_{3} dx_{4} \end{equation}
and the corresponding Euclidean structure is the standard one.  The subspace
$\bR^{4}\subset E_{0}$ is co-associative and any co-associative triple
$(E,\phi,V)$ is equivalent to the model $(E_{0}, \phi_{0}, \bR^{4})$ by an oriented linear isomorphism from $E$ to $E_{0}$.

Somewhat more generally, let $V$ be an oriented 4-dimensional real vector space. The wedge product gives a quadratic form on  2-forms up to a positive factor. Let $H$ be a 3-dimensional vector space and $\uomega \in H^{*}\otimes \Lambda^{2} V^{*}$. We say that $\uomega$ is {\it hypersymplectic} if, regarded as a linear map from $H$ to $\Lambda^{2}V^{*}$, it is an injection to a maximal positive subspace for the wedge product form. There is an induced orientation of $H$ and, with a suitable sign convention, for any oriented volume form ${\rm vol}_{H}\in \Lambda^{3} H^{*}$ the 3-form
\begin{equation}   \phi = \uomega + {\rm vol}_{H}\end{equation}
in $H^{*}\otimes \Lambda^{2} V^{*}\oplus \Lambda^{3} H^{*}\subset \Lambda^{3}(V\oplus H)^{*}$ is positive. Moreover the general positive 3-form on $V\oplus H$ such that $V$ is co-associative and $H$ is the orthogonal complement of $V$ with respect to the induced Euclidean structure has the shape above. Notice that a hypersymplectic $\uomega$ defines a conformal structure on $V$---the unique structure such that the image of $H$ is the self-dual subspace.

Next we consider an oriented $7$-manifold $M$. As stated in the introduction, a {\it closed $G_{2}$-structure} on $M$ is a closed $3$-form $\phi$ which is positive at each point  and the structure is {\it torsion-free} if in addition $d*_{\phi}\phi=0$. If $\phi$ is any positive $3$-form a {\it co-associative submanifold} is a $4$-dimensional submanifold $X\subset M$ such that each point $p\in X$ the tangent space $TX_{p}$ is a co-associative subspace of $TM_{p}$ with respect to $\phi(p)$ \cite{kn:HL}. The form $\phi$ induces a standard orientation of $X$, for example by declaring that $*_{\phi} \phi$ is positive on $X$.  As background, we recall that for general positive forms the co-associative condition is over-determined, but for closed forms $\phi$ there is an elliptic deformation theory  and compact co-associative submanifolds are stable under small perturbations of the $3$-form \cite{kn:McL}. They are parametrised by a moduli space of dimension $b^{+}(X)$ (the dimension of a maximal positive subspace for the intersection form). If $\phi$ is torsion-free then co-associative submanifolds are examples of Harvey and Lawson's {\it calibrated submanifolds}\cite{kn:HL}, with calibrating form $*_{\phi}\phi$.

We pause here to note  a useful identity.
\begin{lem} If $\phi$ is a closed positive form then for any vector $v$ in $TM$ we have
$$   (v\hook \phi)\wedge d*_{\phi} \phi = (v\hook d\phi)\wedge *_{\phi} \phi . $$
\end{lem}
This can be derived from the diffeomorphism invariance of the Hitchin functional
$H(\phi)={\rm Vol}(M, g_{\phi})$. The first variation under a compactly supported variation is
$$ \delta H= \frac{1}{3} \int_{M} (\delta \phi)\wedge *_{\phi}\phi. $$
If $v$ is a compactly supported vector on $M$ and $\delta \phi$ is the Lie derivative $L_{v}\phi$ then clearly the first variation is zero. Writing $L_{v}\phi = d (v\hook \phi) + v\hook d\phi $ and integrating by parts we get the identity
$$\int_{M}( v\hook d\phi) \wedge *_{\phi} \phi- (v\hook \phi)\wedge d*_{\phi} \phi =0, $$ which implies the result.

\

We now move on to our central topic---fibrations $\pi:M\rightarrow B$ of a $7$-manifold over a $3$-dimensional base with co-associative fibres. Later  we will consider \lq\lq fibrations'' in the algebraic geometers sense, with singular fibres, but in this section we will consider a genuine  locally trivial bundle, with fibre a $4$-manifold $X$. For definiteness we take $X$ to be the smooth $4$-manifold underlying a K3 surface, with the standard orientation.  We start with a  fixed $C^{\infty}$ fibration $\pi:M\rightarrow B$ where $M$ and $B$ are oriented. We want to work with differential forms on the total space and to this end we recall some background.

Let $V\subset TM$ be the tangent bundle along the fibres, so we have an exact sequence
$$  0\rightarrow V\rightarrow TM\rightarrow \pi^{*}(TB)\rightarrow 0.$$ 
For $p,q\geq 0$ let $\Lambda^{p,q}$ be the vector  bundle 
$$\Lambda^{p}(\pi^{*} T^{*}B)\otimes \Lambda^{q} V^{*}$$
over $M$ 
and write $\Omega^{p,q}$ for the space of sections of $\Lambda^{p,q}$. The exterior derivative on the fibres  is a natural differential operator, for $p,q\geq 0$:
$$d_{f}: \Omega^{p,q}\rightarrow \Omega^{p,q+1},$$ and there is a canonical filtration of $\Omega^{*}_{M}$ with quotients the $\Omega^{*,*}$. 

Now suppose that we have a {\it connection} $H$ on $\pi:M\rightarrow B$. That is, a sub-bundle $H\subset TM$ with $TM=V\oplus H$. Then we get a direct sum decomposition
$$  \Omega^{*}_{M}= \bigoplus_{p,q\geq 0} \Omega^{p,q}. $$
The exterior derivative $d:\Omega^{*}_{M}\rightarrow \Omega^{*}_{M}$ has three components with respect to this decomposition.
\begin{enumerate} \item The exterior derivative along the fibres $d_{f}$, as above. 
   \item The coupled exterior derivative in the horizontal direction
   $$ d_{H}: \Omega^{p,q}\rightarrow \Omega^{p+1,q}. $$
This can be defined in various ways. One way is to  work at a given point $b\in B$ and to choose a local trivialisation of the fibration compatible with $H$ on $\pi^{-1}(b)$. Then, over this fibre, $d_{H}$ is given  by the horizontal component of the exterior derivative in this trivialisation. From another point of view we can regard $H$ as a connection on a principal bundle with structure group ${\rm Diff}(X)$. Then $d_{H}$ is  the usual coupled exterior derivative on an infinite dimensional associated  vector bundle.   
\item The curvature term. The curvature of the connection $H$ is a section $F_{H}$ of the bundle  $V\otimes \pi^{*}\Lambda^{2} T^{*} B$ defined by the equation
$$    [\xi_{1}, \xi_{2}]_{V}= F_{H}(\xi_{1}, \xi_{2}), $$
for any sections $\xi_{1}, \xi_{2}$ of $H$, where $[\ ,\ ]_{V}$ denotes the vertical component of the Lie bracket. (This agrees with ordinary notion of curvature if we take the point of view of a ${\rm Diff}(X)$ bundle.) The tensor $F_{H}$ acts algebraically on $\Omega^{p,q}$ by wedge product on the horizontal term and contraction on the vertical term,  yielding a map (which we denote by the same symbol)
$$  F_{H}: \Lambda^{p,q}\rightarrow \Lambda^{p+2,q-1}. $$

\end{enumerate}
The fundamental differential-geometric formula for the exterior derivative on the total space is then
\begin{equation}   d= d_{f}+d_{H}+ F_{H}: \Omega^{p,q}\rightarrow \Omega^{p,q+1} + \Omega^{p+1,q} +\Omega^{p+2, q-1}. \end{equation}

   We say that $\uomega\in \Omega^{1,2}$ is a hypersymplectic element if it is hypersymplectic at each point, as defined above. 
\begin{prop}
Given $\pi:M\rightarrow B$ as above, a  closed $G_{2}$ structure on $M$ such that the fibres of $\pi$ are co-associative, with orientation compatible with those given on $M $and  $B$, is equivalent to the following data. 
\begin{enumerate}
\item A connection $H$.
\item A hypersymplectic element  $\uomega\in \Omega^{1,2}$ with $d_{f}\uomega=0, d_{H}\uomega=0$.        
\item A tensor $\ulambda \in \Omega^{3,0}$ with 
\begin{equation} d_{f} \ulambda = -F_{H} \uomega \end{equation}
and such that the value of $\ulambda$ at each point is positive, regarded as an element of $\Lambda^{3}T^{*}B$. 
\end{enumerate}

\end{prop}
 
 This proposition follows immediately from the algebraic discussion at the beginning of this section and the description of the exterior derivative above. The corresponding  positive $3$-form $\phi$ is of course $\uomega+\ulambda$.

Note that if we have {\it any} $3$-form $\phi$ which vanishes on the fibres there is a well-defined component $\uomega$ of $\phi$ in $\Omega^{1,2}$, independent of a connection.  This just reflects the naturality of the filtration
of $\Omega^{*}_{M}$. If $\phi$ is closed then $d_{f}\uomega=0$. At this point it is convenient to switch to a local discussion (in the base), so we suppose for now that $B$ is a $3$-ball with co-ordinates $t_{1}, t_{2}, t_{3}$. Suppose we have any closed $3$-form $\phi$ on $M$ which vanishes on the fibres. Then the de Rham cohomology class of $\phi$ is zero and we can write $\phi= d\eta$ for a $2$-form $\eta$ on $M$. By assumption, the restriction of $d \eta$ to each fibre is zero so we get a map
\begin{equation} h: B\rightarrow H^{2}(X), \end{equation} taking a point to the cohomology class of $\eta$ in the fibre over that point.
Changing the choice of $\eta$ changes $h$ by the  addition of a constant in $H^{2}(X)$.  As discussed above, the form $\phi$ has a well-defined component $\uomega \in \Omega^{1,2}$ which we can write
$$  \uomega= \sum_{i=1}^{3} \omega_{i} dt_{i}, $$
 where the $\omega_{i}$ are closed $2$-forms on the fibres. If we fix a trivialisation of the bundle then the $\omega_{i}$ become families of closed $2$-forms on a fixed 4-manifold $X$, parametrised by $(t_{1}, t_{2}, t_{3})$. A moments thought shows that the de Rham cohomology class $[\omega_{i}]$ is 
\begin{equation} [\omega_{i}]= \frac{\partial h}{\partial t_{i}}. \end{equation}
This formula makes intrinsic sense, since the cohomology groups of the different fibres are canonically identified. We say that a map $h:B\rightarrow H^{2}(X,\bR)$ is {\it positive} if it is an immersion and the image of its derivative at each point is a positive subspace with respect to the cup product form. Then it is clear that a necessary condition for the existence of a hypersymplectic element $\uomega$ representing the derivative is that $h$ is positive.

\

\

We now change point of view and ask how we can build  up a form $\phi$.
We start  with a positive map $h:B\rightarrow H^{2}(X)$ and suppose that we have chosen a hypersymplectic element  $\uomega=\sum \omega_{i} dt_{i}\in \Omega^{1,2}$ representing the derivative of $h$.  
\begin{lem}
Given $\uomega$ as above there is a connection $H$ such that $d_{H}\uomega=0$.
\end{lem}

We work in a local trivialisation so that the $\omega_{i}$ are regarded as $\ut$-dependent $2$-forms on the fixed $4$-manifold $X$. The connection is represented by $\ut$-dependent vector fields  $v_{1}, v_{2}, v_{3}$ on $X$ and the condition that $d_{H}\uomega=0$ is
\begin{equation}  \frac{\partial \omega_{i}}{\partial t_{j}} - \frac{\partial \omega_{j}}{\partial t_{i}} + (L_{v_{j}}\omega_{i}-  L_{v_{i}}\omega_{j}) =0 \end{equation} 
By assumption, $\omega_{i}$ represents the cohomology class $\frac{\partial h}{\partial t_{i}}$ and this implies that the cohomology class of  $\frac{\partial \omega_{i}}{\partial t_{j}} - \frac{\partial
\omega_{j}}{\partial t_{i}}$ is zero. So we can choose $\ut$-dependent  $1$-forms $a_{1}, a_{2}, a_{3}$ such that
$$   \frac{\partial \omega_{i}}{\partial t_{j}} - \frac{\partial
\omega_{j}}{\partial t_{i}}= da_{k}, $$
for $i,j,k$ cyclic. Since the $\omega_{i}$ are closed we can write (9) as
$$  d( v_{j}\hook \omega_{i}- v_{i}\hook \omega_{j}) =da_{k}, $$
which is certainly satisfied if
\begin{equation}  v_{j}\hook \omega_{i}- v_{i}\hook \omega_{j}= a_{k} \end{equation}
Let $$ S: V\oplus V\oplus V\rightarrow V^{*}\oplus V^{*}\oplus V^{*} $$
 be the map which takes a triple of tangent vectors $v_{1}, v_{2}, v_{3}$ to a triple of co-tangent vectors $v_{j}\hook \omega_{i}-v_{i}\hook \omega_{j}$. What we need is the elementary linear algebra statement that $S$ is an isomorphism.
To see this observe that it is  a pointwise statement and by making a linear change of co-ordinates $t_{i}$ we can reduce to the case when $\omega_{i}$ are the standard basis of the self-dual forms on $\bR^{4}$. Using the metric to identify tangent vectors and cotangent vectors,  $S$ becomes the map which takes a triple $v_{1}, v_{2}, v_{3}$ to the triple \begin{equation}    J v_{3}- K v_{2}\ ,\  K v_{1}- I  v_{3} \ ,\  I v_{2}- J v_{1}. \end{equation}
Here we make the standard identification of $\bR^{4}$ with the quaternions and use quaternion multiplication. 
A short calculation in quaternion linear algebra shows that  $S$ is an isomorphism.

Once we have chosen $\uomega$ and $H$, as above,  the remaining data we need to construct a closed $G_{2}$ structure with co-associative fibres on $M$ is $\ulambda\in \Omega^{3,0}$ which we can write in co-ordinates as $-\lambda dt_{1} dt_{2} dt_{3}$ for a function $\lambda$ on $M$. The identity $d^{2}=0$ on $\Omega^{*}_{M}$ implies that
$$   d_{f} ( F_{H} \uomega) =0. $$
Since $H^{1}(X)=0$ we can solve the equation
$$   d_{f} \ulambda= - F_{H}\uomega. $$
The solution is unique up to the addition of a lift of a $3$-form  from $B$, i.e. up to changing $\lambda$ by the addition of a function of $\ut$. By adding a sufficiently positive function of $\ut$ we can arrange that
$\lambda$ is positive. We conclude that  any choice of hypersymplectic element representing the derivative of $h$ can be extended to a closed $G_{2}$ structure with co-associative fibres.

\

We stay with the local discussion and go on to consider torsion-free $G_{2}$ structures. Suppose we have $2$-forms $\omega_{i}$ spanning a maximal positive subspace of the vertical tangent space $V$. Let $\chi$ be an arbitrary volume form on $V$, so that $(\omega_{i}\wedge\omega_{j})/\chi$ are real numbers. We form the determinant $\det((\omega_{i}\wedge \omega_{j})/\chi)$ and the $4$-form $$ \left({\rm det} ((\omega_{i}\wedge\omega_{j})/\chi))\right)^{1/3} \chi. $$ It is clear that this is independent of the choice of $\chi$ and we just write it as ${\rm \det}^{1/3}(\omega_{i}\wedge \omega_{j})$. 
\begin{lem}
If $\phi= \uomega+ \ulambda \in  \Omega^{1,2}+ \Omega^{3,0}$ is a positive $3$-form on $M$, as considered above,  then
$$*_{\phi} \phi= \uTheta+ \umu$$ where
\begin{enumerate}
\item $\uTheta={\cal F}_{1}(\uomega,\ulambda)\in \Omega^{2,2}$ is given by 
$ \sum_{\rm cyclic} \Theta_{i} dt_{j}dt_{k} $ where $\Theta_{i}$ are forms in the linear span of $\omega_{1},\omega_{2}, \omega_{3}$ determined uniquely by the condition that \begin{equation}
\Theta_{i} \wedge \omega_{j} = \delta_{ij} (\lambda^{1/3} {\rm det}^{1/3}(\omega_{a}\wedge\omega_{b}));\end{equation}
\item $\umu= {\cal F}_{2}(\uomega,\ulambda)\in \Omega^{4,0}$ is given by
\begin{equation}   {\rm det}^{1/3}(\omega_{a}\wedge\omega_{b})\lambda^{-2/3}. \end{equation}

\end{enumerate}\end{lem}

This is a straightforward algebraic exercise using the definitions. The upshot of our local discussion is then the following statement. 
\begin{prop}
 A torsion-free $G_{2}$ structure on $M$ with co-associative fibres is given by a map $h$ (defined up to a constant), a  hypersymplectic element $\uomega$ representing the derivative of $h$, a connection  $H$ and  $\ulambda\in \Omega^{3,0}$ with $ \ulambda>0$, satisfying the five equations
\begin{enumerate}  
\item $d_{H} \uomega=0$,
\item $d_{f}\ulambda= -F_{H} \uomega$,
\item $d_{H} \umu =0$,
\item $d_{f} \uTheta =- F_{H} \umu$, 
\item $d_{H}\uTheta=0$.
\end{enumerate}
Here $\umu, \uTheta$ are determined algebraically from $\uomega, \ulambda$ by the formulae (12),(13). 
\end{prop}

Note that the equation (3) in Proposition 2 has a clear geometric meaning. It states that the structure group of the connection reduces to the volume preserving diffeomorphisms of $X$. It is also equivalent to the statement that the fibres are minimal submanifolds, with respect to the metric $g_{\phi}$.

We will not attempt to analyse these equations further in this form, except to note one identity.
\begin{lem}
Suppose that $h,\uomega, \ulambda, H$ satisfy equations (1), (2), (3), (4) in Proposition 2. Then $d_{H}\uTheta= G dt_{1} dt_{2} dt_{3}$ where $G$ is an anti-self dual harmonic form on each fibre with respect to the conformal structure determined by $\uomega$. 
\end{lem}
 First, the identity $d^{2}*_{\phi}\phi=0$ implies that, under the stated conditions $d_{f}(d_{H}\uTheta)=0$ so $d_{H}\Theta=G dt_{1}dt_{2} dt_{3}$ where $G$ is a closed $2$-form on each fibre. Now we apply Lemma 1,  taking $v$ to be a horizontal vector. Since $d\phi=0$ we have $(v\hook\phi)\wedge d*_{\phi}\phi=0$ . The contraction $v\hook\phi$ has terms in $\Omega^{0,2}$ and $\Omega^{2,0}$ but only the first contributes to the wedge product with $d *_{\phi} \phi$ and the result follows immediately.  

\

Let us now move back to consider the \lq\lq global'' case, with a K3 fibration $\pi:M\rightarrow B$ over a general oriented $3$-manifold $B$. We have a flat bundle ${\cal H}\rightarrow B$ with fibre $H^{2}(X;\bR)$ induced from an $H^{2}(X;\bZ)$ local system, so there is an intrinsic integer lattice in each fibre.  Let  $\phi$ be a closed $3$-form on $M$, vanishing on the fibres. If we cover $B$ by co-ordinate balls $B_{\alpha}$ we can analyse  $\phi$  on each $\pi^{-1}(B_{\alpha})$ as above, starting with a choice of a smooth section $h_{\alpha}$ of ${\cal H}\vert_{B_{\alpha}}$. Then on $B_{\alpha}\cap B_{\beta}$ the difference 
\begin{equation}\chi_{\alpha \beta}= h_{\alpha}-h_{\beta}\end{equation}  is a constant section of ${\cal H}$. This defines a Cech cocycle, and hence a cohomology class
$$  \chi \in H^{1}(B; \underline{{\cal H}}) , $$
where $\underline{{\cal H}}$ is the sheaf of locally constant sections of the flat bundle ${\cal H}$.  
This class $\chi$ is determined by the de Rham cohomology class of $\phi$. In fact the Leray-Serre spectral sequence gives us an exact sequence
   \begin{equation}    0\rightarrow H^{3}(B)\rightarrow H^{3}(M)\stackrel{b}{\rightarrow} H^{1}(B; \underline{{\cal H}})\rightarrow 0 ,  \end{equation}
and $[\phi]\in H^{3}(M)$ maps to $\chi\in H^{1}(B;\underline{{\cal H}})$. It is convenient to express this in terms of an \lq\lq affine'' variant of ${\cal H}$. Let $G$ be the group of affine isometries of $H^{2}(X)$, an extension
$$ 1 \rightarrow    H^{2}(X)\rightarrow G \rightarrow O(H^{2}(X))\rightarrow 1. $$
A class $\chi\in H^{1}(B; \underline{{\cal H}})$ defines a flat bundle ${\cal H}_{\chi}$ with fibre $H^{2}(X)$ and  structure group $G$ over $B$. One way of defining this is to fix a cover and Cech representative $\chi_{\alpha \beta}$ and then decree that locally constant sections of ${\cal H}_{\chi}$ are given by constant sections $h_{\alpha}$ over $B_{\alpha}$ satisfying (14). In these terms the conclusion is that if we start with a class $\Phi\in H^{3}(M)$,  map to a class $\chi=b(\Phi)$ and form the bundle ${\cal H}_{\chi}$, then any representative $\phi\in \Omega^{3}(M)$ for $\Phi$ which vanishes on the fibres defines a canonical smooth section $h$ of the flat bundle ${\cal H}_{\chi}$. Moreover if $\phi$ is a positive $3$-form then $h$ is a {\it positive} section of ${\cal H}_{\xi}$, using the obvious extension of the definition in the local situation.  

\subsection{The adiabatic limit}

 We have seen that the torsion-free condition for a positive $3$-form involves six equations for a triple $(\uomega,\ulambda, H)$
$$\begin{array}{lll}
 d_{f}\uomega=0, & 
  d_{H} \uomega=0,&
 d_{f}\ulambda= -F_{H} \uomega,\\
d_{H} \umu =0,& d_{f} \uTheta = -F_{H} \umu, 
& d_{H}\uTheta=0.
\end{array}$$
where $\uTheta={\cal F}_{1}(\uomega,\ulambda)$ and $\umu={\cal F}_{2}(\uomega,\ulambda)$ are determined algebraically by $\uomega,\ulambda$.  
We introduce a positive real parameter $\epsilon$ and write $\uomega=\epsilon\tuomega$,
$$\tuTheta= {\cal F}_{1}(\tuomega, \lambda)= \epsilon\uTheta, \ \ \  \tumu={\cal F}_{2}(\tuomega,\lambda)=\epsilon^{2}\umu. $$ 
 Then the six equations, written for $(\tuomega, \ulambda, H)$ become
\begin{equation} \begin{array}{lll}
 d_{f}\tuomega=0, & 
  d_{H} \tuomega=0,&
 d_{f}\ulambda+\epsilon F_{H} \tuomega=0,\\
d_{H} \tumu =0,& d_{f} \tuTheta  +\epsilon F_{H} \tumu=0, 
& d_{H}\tuTheta=0.
\end{array}\end{equation}

Here $\tumu,\tuTheta$ are related to $\tuomega,\ulambda$ in just the same way as $\umu,\uTheta$ were to $\uomega,\ulambda$. Geometrically, in the metric determined by $\tuomega,\ulambda$ the volume of the fibre is $\epsilon^{2}$ times that in the metric determined by $\uomega,\ulambda$. For non-zero $\epsilon$ these equations are   completely equivalent to the original set; the fundamental idea we want to pursue is to pass to the \lq\lq adiabatic limit'', setting  $\epsilon=0$ so that the terms involving the curvature $F_{H}$ drop out.  At this stage we simplify notation by dropping the tildes, writing $\uomega,\uTheta,\umu$ for $\tuomega,\tuTheta,\tumu$.

\

We say that  $\uomega\in \Omega^{1,2}$ is a {\it hyperk\"ahler element} if, in local co-ordinates on the base, $\uomega=\sum \omega_{i}dt_{i}$ with
$$\omega_{i}\wedge \omega_{j} = a_{ij} \nu, $$
where $\nu$ is a volume form along the fibres and $a_{ij}$ is a positive definite matrix {\it constant} along the fibres. 
\begin{lem}
Suppose that a pair $(\uomega,\ulambda)$ satisfies the five equations
$$ \begin{array}{lll}
 d_{f}\uomega=0, & 
  d_{H} \uomega=0,&
 d_{f}\ulambda= 0 ,\\
d_{H} \umu =0,& d_{f} \uTheta = 0. 
& 
\end{array}$$

Then $\uomega$ is a hyperk\"ahler element and $\ulambda$ is the lift of a positive $3$-form on $B$. Conversely given a hyperk\"ahler element $\uomega$ and a positive $3$-form on $B$ there is a unique connection $H$ which yields a solution to these equations.
\end{lem}

In one direction, suppose that we have $\uomega, \ulambda$ satisfying these equations and work locally,  writing $\uomega=\sum \omega_{i} dt_{i}, \ulambda=-\lambda dt_{1} dt_{2} dt_{3}, \uTheta=\sum \Theta_{i} dt_{j} dt_{k}$. The equation $d_{f}\ulambda=0$ immediately gives that $\ulambda$ is a lift from $B$.  We know that $\uTheta_{i}$ takes values in the self-dual space  $\Lambda^{+}$ spanned by  $\omega_{1}, \omega_{2}, \omega_{3}$ so the equation $d_{f}\uTheta=0$ says that the $\Theta_{i}$ are self-dual harmonic forms. Since $b^{+}(X)=3$ we have $\Theta_{i}=\sum g_{ij} \omega_{j}$ where the functions $g_{ij}$ are constant on the fibres. Now
$$   \Theta_{i}\wedge \omega_{k}= \lambda^{1/3}\ \delta_{ik}\ {\rm det}^{1/3}(\omega_{a}\wedge \omega_{b}), $$
so $$  \sum_{j} g_{ij} \omega_{j}\wedge \omega_{k}= \lambda^{1/3}\ \delta_{ik}
\ {\rm det}^{1/3}(\omega_{a}\wedge
\omega_{b}) $$ which implies that $$\omega_{j}\wedge \omega_{k} = g^{jk} \lambda^{1/3} {\rm det}^{1/3}(\omega_{a}\wedge
\omega_{b})$$
where $g^{jk}$ is the inverse matrix, also constant on the fibres. Hence $\uomega$ is hyperk\"ahler.

For the converse, suppose that we have a hyperk\"ahler element $\uomega$ and volume form on $B$, which we may as well take to be $-dt_{1}dt_{2} dt_{3}$ in local co-ordinates. Then we have a volume form $\umu$ along the fibres given by ${\rm det}^{1/3}(\omega_{i}\wedge \omega_{j})$ and we have to show that there is a unique connection $H$ such that  $d_{H}\uomega=0$ and $d_{H}\umu=0$. For this we go back to the proof of Lemma 2. We can choose a local trivialisation such that the volume form along the fibres is constant. For any choice of $1$-forms $a_{1}, a_{2}, a_{3}$ with 
$$ \frac{\partial \omega_{i}}{\partial t_{j}}- \frac{\partial \omega_{j}}{\partial t_{i}}= da_{k}$$ we showed there that there is a unique choice of connection $H$. But we are free to change the $1$-forms $a_{k}$ by the derivatives of functions. It is easy to check that the condition that the connection is volume preserving is exactly that the $1$-forms $a_{k}$ satisfy $d^{*}a_{k}=0$, with respect to the induced metric on the fibres. By Hodge theory,  this gives the existence and uniqueness of the connection $H$.  (In fact the argument   shows that for any volume form along the fibres we can choose $H$ to satisfy this additional condition, provided that the volume of the fibres is constant.)

\

We can clearly choose the solution in Lemma 5 so that the $\umu$-volume of the fibres is $1$, and we fix this normalisation from now on.

\

It remains to examine the last equation, $d_{H}\uTheta=0$. Recall that 
$H^{2}(X,\bR)$ carries an intrinsic quadratic from of signature $(3,19)$. This defines a volume form on any submanifold and hence an Euler-Lagrange equation defining stationary submanifolds. In our context we are interested in $3$-dimensional stationary submanifolds whose tangent spaces are  positive subspaces
and these are called  {\it maximal submanifolds}. A section $h$ of ${\cal H}_{\xi}$ can be represented locally by an embedding of the base in $H^{2}(X)$ and we say that $h$ is a maximal section if the image is a maximal submanifold. More generally, for any section $h$ the {\it mean curvature} $m(h)$ is defined as a section of the vector bundle $({\rm Im}\  dh)^{\perp} \subset {\cal H}$ which  is  the normal bundle of the submanifold. 

\begin{lem}
Suppose that $h$ is a positive section of ${\cal H}_{\chi}$,  $\uomega$ is a hyperk\"ahler element representing the derivative of $h$, and that $H$, $\ulambda$ are chosen as above. 
Then \begin{equation} d_{H}\uTheta= m(h) \ulambda,  \end{equation} where we use the metric on the fibres defined by $\uomega$ to identify $({\rm Im}\ dh)^{\perp}$ with the anti-self-dual 2-forms along the fibres and use the result of Lemma 4. In particular, $d_{H}\uTheta=0$ if and only if $h$ is a maximal positive section. 
\end{lem}

In local co-ordinates $t_{i}$ let $g_{ij}= \int_{X}\omega_{i}\wedge \omega_{j}$ and set $J=\det(g_{ij})$. The volume of the image is given by
$$  \int J^{1/2} dt_{1} dt_{2} dt_{3} $$ and the mean curvature is
$$     m= J^{-1/2} \left( J^{1/2} g^{ij} h_{,j}\right)_{,i},  $$
where $g^{ij}$ is the inverse matrix and $\ _{,i}$ denote partial differentiation. That is
$$    m =  J^{-1/2} \sum \frac{\partial \Theta^{*}_{i}}{\partial t_{i} }, $$
where $\Theta^{*}_{i}= J^{1/2} g^{ij} [\omega_{j}]$. Computing in a trivialisation compatible with $H$ at a given point we have
$$    d_{H}\uTheta= \left(\sum \frac{\partial \Theta_{i}}{\partial t_{i} }\right) dt_{1}dt_{2} dt_{3}. $$
Now $\Theta_{i}$ is defined by the identity 
$$   \Theta_{i}\wedge \omega_{j}= {\rm det}^{1/3} \lambda^{1/3} \delta_{ij},$$ so

$$   \int_{X} \Theta_{i}\wedge \omega_{j}= \lambda^{1/3} \int_{X} {\rm det}^{1/3} \delta_{ij}, $$
 whereas $$ \int_{X} \Theta^{*}_{i}\wedge \omega_{j}= J^{1/2}. $$
 The condition that the $\mu$-volume of the fibres is $1$ is 
 $$  \lambda^{2/3} =\int_{X} {\rm det}^{1/3}. $$ On the other hand one has
 $$    \int_{X} {\rm det}^{1/3} = J^{1/3}. $$
 This implies that $\Theta^{*}_{i}=[\Theta_{i}]$ and $J^{-1/2} dt_{1} dt_{2} dt_{3} = \ulambda$  and the result follows.

\

\

One can also prove this Lemma more synthetically, comparing the Hitchin 7-dimensional volume functional with the 3-dimensional volume of the image of $h$.

\subsection{Formal power series solutions} 

The conclusion that we are lead to  in the preceding discussion is that the adiabiatic limit of the torsion-free $G_{2}$-structure condition should be the maximal submanifold equation $m(h)=0$.  To give more support to this idea we show that a solution  of the latter equation can be a extended to a formal power series solution of system of equations (16). For simplicity we work locally, so we take the base $B$ to be a ball.
\begin{prop}\begin{enumerate} \item Suppose that the quadruple $(h,\uomega, H, \ulambda)$ satisfies the five equations of Lemma 5.  Then there are formal power series
$$\tuomega^{\epsilon}=\uomega+\sum_{j=1}^{\infty} \uomega^{(j)}\epsilon^{j},\  H^{\epsilon}= H+\sum_{j=1}^{\infty} \zeta^{(j)}\epsilon^{j}, \ \ulambda^{\epsilon}=\ulambda+\sum_{j=1}^{\infty} \ulambda^{(j)}\epsilon^{j}$$
which give a formal power series solution of  the first five equations of (16). We can choose $\uomega^{(j)}$ to be exact, so that $\tuomega^{\epsilon}$ represents the derivative of the fixed map $h$.
\item Suppose that in addition $h$ satisfies the maximal submanifold equation $m(h)=0$ and is smooth up to the boundary of the ball $B$.  Then we can find a formal power series 
$$  h^{\epsilon}= h+\sum_{j=1}^{\infty} h^{(j)} \epsilon^{j} $$
 and power series $\tuomega^{\epsilon}, H^{\epsilon}, \umu^{\epsilon}$ as above but now with $\tuomega^{\epsilon}$ representing the derivative of $h^{\epsilon}$ and such that the quadruple $(h^{\epsilon},\tuomega^{\epsilon},
H^{\epsilon}, \ulambda^{\epsilon})$ is a formal power series solution of the equations (16). \end{enumerate}
\end{prop}

The proof follows standard lines once we understand the linearisation of the equations, and for this we to develop some background. Let $\omega_{i}$ be a standard hyperk\"ahler triple on $X$. Consider an infinitesimal variation $\omega_{i}+ \eta_{i}$ with $d\eta_{i}=0$ and write $\eta_{i}= \eta^{+}_{i}+\eta_{i}^{-}$ in self-dual and anti-self-dual parts. So $\eta^{+}_{i}=\sum S_{ij} \omega_{j}$ for some co-efficients $S_{ij}$. A short calculation shows that the variation in $d_{f}\uTheta$ is, to first order, $\sum dT_{ij}\wedge \omega_{j}$  where
\begin{equation} T_{ij} =  S_{ij}+S_{ji}- \frac{2}{3} {\rm Tr}(S) \delta_{ij}. \end{equation}
(i.e. $T$ is obtained from $S$ by projecting to the symmetric, trace-free component.) We want to vary $\omega_{i}$ in fixed cohomology classes, by $da_{i}$. The resulting first order variation in $d\Theta_{i}$ is given by a linear map ${\cal L}_{1}$ from $\Omega^{1}_{X}\otimes \bR^{3}$ to $\Omega^{3}_{X}\otimes \bR^{3}$. So the calculation above shows that this linear map is a composite
\begin{equation} \Omega^{1}_{X}\otimes \bR^{3}\stackrel{d^{+}\otimes 1}{\rightarrow}\Omega^{2,+}_{X}\otimes \bR^{3} \stackrel{A}{\rightarrow} \Omega^{2,+}_{X}\otimes \bR^{3}\stackrel{d\otimes 1}{\rightarrow}\Omega^{3}_{X}\otimes \bR^{3},  \end{equation}
where $A$ is the map taking $\sum S_{ij} \omega_{j}$ to $\sum T_{ij}\omega_{j}$, with $T$ given by (18) above. Thus the image of $A$ is the subspace of $\Omega^{2,+}_{X}\otimes \bR^{3}$ consisting of  tensors with values in $s^{2}_{0}(\bR^{3})$---symmetric and trace-free. Here of course we are using repeatedly the fact that the bundle of self-dual forms is identified with $\bR^{3}$ by the sections $\omega_{i}$.

From another point of view, write $S_{+}, S_{-}$ for the spin bundles over $X$ and write $S^{p}_{+}, S^{p}_{-}$ for their pth. symmetric powers.  Of course in our situation $S_{+}$ is a trivial flat bundle. The tensors with values in $s^{2}_{0}(\bR^{3})$ can be identified with sections of $S^{4}_{+}$.  (More precisely, $S^{4}_{+}$ has a real structure and we should consider real sections, but in the discussion below we do not need to distinguish between the real and complex cases.) There is a contraction map from $\Omega^{3}_{X}\otimes \bR^{3}$ to $\Omega^{3}_{X}$ given by
\begin{equation} P(\psi_{1}, \psi_{2}, \psi_{3})= \sum (*\psi_{i})\wedge \omega_{i} \end{equation}
The kernel of this contraction map can be identified with $S^{3}_{+}\otimes S_{-}$, corresponding to the decomposition:
$$  (S_{+}\otimes S_{-})\otimes S_{+}^{2}= (S_{+}\otimes S_{-})\oplus (S^{3}_{+}\otimes S_{-}). $$
Starting with $$d\otimes 1: \Omega^{2,+}_{X}\otimes \bR^{3}\rightarrow \Omega^{3}_{X}\otimes \bR^{3}, $$
we first restrict to sections of $S_{+}^{4}$ and then project to sections of $S_{+}^{3}\otimes S_{-}$. This defines a differential operator
$$  \delta_{1}:\Gamma (S^{4}_{+})\rightarrow \Gamma(S^{3}_{+}\otimes S_{-}). $$
On the other hand, starting with
$$  d:\Omega^{3}_{X}\otimes \bR^{3}\rightarrow \Omega^{4}_{X}\otimes \bR^{3}, $$
we restrict to sections of $S^{3}_{+}\otimes S_{-}$ to define
$$   \delta_{2}:\Gamma(S^{3}_{+}\otimes S_{-}) \rightarrow \Gamma(S_{+}^{2}). $$
So we have $$  \Gamma(S^{4}_{+})\stackrel{\delta_{1}}{\rightarrow} \Gamma(S^{3}_{+}\otimes S_{-})\stackrel{\delta_{2}}{\rightarrow}\Gamma(S^{2}_{+}). $$
The key fact we need is that this is exact at the middle term, i.e. that
\begin{equation} {\rm Ker} \ \delta_{2}={\rm Im} \ \delta_{1}. \end{equation}
In fact, calculation shows that $\delta_{2}+ \delta_{1}^{*}$ can be identified with the Dirac operator coupled with $S^{3}_{+}$,  i.e.

$$   D:\Gamma (S^{3}_{+}\otimes S_{-})\rightarrow \Gamma(S^{3}_{+}\otimes S_{+})$$
using the isomorphism
$$   S_{+}^{3}\otimes S_{+}= S_{4}^{+}\oplus S_{2}^{+}. $$
Since $S_{+}$ is trivial this Dirac operator is a sum of copies of the ordinary Dirac operator over the K3 surface and has no kernel in the $S_{-}$ term. This implies (21). We can say a bit more. The kernel of $\delta_{1}$ consists of the constant sections of $S_{+}^{4}$ so for any $\rho\in S_{3}^{+}\otimes S_{-}$ we can solve the equation $\delta_{1} T=\rho$ with $T$ in the  $L^{2}$-orthogonal complement of the constants. But this means that we can solve the equation $(d^{+}\otimes 1) a = T$.

 To sum up we have shown:
 \begin{lem}If  $\rho\in \Omega_{X}^{3}\otimes\bR^{3}$ satisfies the two conditions
 \begin{itemize}\item $P(\rho)=0$ where $P$ is the contraction (20);
 \item ($d\otimes 1)(\rho)=0$ in $\Omega_{X}^{4}\otimes \bR^{3}$; \end{itemize}then $\rho$ is in the image of ${\cal L}_{1}:  \Omega_{X}^{1}\otimes \bR^{3}\rightarrow \Omega_{X}^{3}\otimes \bR^{3}$.

\end{lem}

\

 We now turn to the proof of the first item in Proposition 3. At stage $k$ we suppose that we have found $\uomega^{(i)}, \ulambda^{(i)}, \zeta^{(i)}$ for
$i\leq k$ so that  the finite sums
$$  \uomega^{[k]}=\uomega+\sum_{i=1}^{k} \uomega^{(k)} \ \ \ {\rm etc.}, $$
satisfy the first five equations of (16) modulo $\epsilon^{k+1}$. We want to choose $\uomega^{(k+1)}, \ulambda^{(k+1)} $ and $ \zeta^{(k+1)}$ so that if we define
$$  \uomega^{[k+1]}= \uomega^{[k]}+ \epsilon^{k+1} \uomega^{(k+1)} \ \ \ {\rm etc.}$$
we get solutions of the first five equations of (16) modulo $\epsilon^{k+2}$. The essential point is that this involves solving {\it linear} equations for $\uomega^{(k+1)}, \ulambda^{(k+1)}, \zeta^{(k+1)}$, and moreover these equations are essentially defined fibrewise.

\

{\bf Step 1}

Let $E_{3,1}$ be the \lq\lq error term'' 
$$   d_{f}\ulambda^{[k]}+\epsilon F_{H^[k]} \uomega^{[k]} = E_{3,1}\epsilon^{k+1} +O(\epsilon^{k+2}).$$
(Here the $O(\ )$ notation means in the sense of formal power series expansion.)
The equation we need to solve to correct this error term   is
$$    d_{f} \ulambda^{(k+1)}= - E_{3,1},  $$
involving only $\ulambda^{(k+1)}$. 
The identity $d^{2}=0$ on $M$ implies that
$$   d_{f}   F_{H^[k]} \uomega^{[k]}= -d_{H}d_{H^{[k]}}\uomega^{[k]}$$
which is $O(\epsilon^{k+1})$. This means that $d_{f}E_{3,1}=0$ so we can solve the equation for $\ulambda^{(k+1)}$. The solution is unique up to the addition of a term which is constant on the fibres. 

\

{\bf Step 2}

Now let $E_{2,3}$ be the error term
$$  d_{f}\uTheta^{[k]}+\epsilon  F_{H^[k]} \umu^{[k]}= E_{2,3} \epsilon^{k+1}+O(\epsilon^{k+2}). $$
The equation we need to solve to correct this error involves the first order variation of $d_{f}\uTheta^{[k]}$ with respect to changes in $\uomega$ and $\ulambda$. For the first, we write $\omega^{(k+1)}$ as the derivative of $a\in \Omega^{1}\otimes \bR^{3}$ as above and the first order variation is given by ${\cal L}_{1}(a)$. Denote the variation with respect to $\ulambda$ by ${\cal L}_{2}$. The equation we need to solve is
$$   {\cal L}_{1}(a)+ {\cal L}_{2}(\ulambda^{(k+1)}) = -E_{2,3}. $$
The same argument as before shows that $d_{f}E_{2,3}=0$. We do not need to write our ${\cal L}_{2}$ explicitly because it is easy to see that ${\cal L}_{2}(\ulambda^{(k+1)})$ only involves the derivative along the fibre, so is unchanged if we change $\ulambda^{(k+1)}$ by a term which is constant on the fibres. Thus  ${\cal L}_{2}(\ulambda^{(k+1)})$ is fixed by Step 1 and we have to solve the equation
$$    {\cal L}_{1}(a)= \rho, $$
for $a$, 
where $ \rho=-E_{2,3}- {\cal L}_{2}(\ulambda^{(k+1)})$.
 The same argument as before shows that $d_{f}E_{2,3}=0$ and it is clear that $d_{f} \circ {\cal L}_{2}=0$. Thus $d_{f}\rho=0$. We claim that, on each fibre, $P(\rho)=0$. Given this claim, Lemma 7 implies that we can find the desired solution $a$.

To verify the claim we use the identity of Lemma 1. We take $v$ to be a vertical vector and $\phi$ to be the $3$-form
$$ \phi= \omega^{[k]} + \ulambda^{[k+1]} . $$
By what we have arranged, the $(3,1)$ component of $d\phi$ is $O(\epsilon^{k+2})$ and this means that 
$  (v\hook \uomega^{[k]})\wedge \rho $ is $O(\epsilon)$. In turn this implies that  $(v\hook \uomega)\wedge \rho$ vanishes for all vertical vectors $v$,  and it is easy to check that this is equivalent to the condition $P(\rho)=0$.

\

{\bf Step 3} At this stage we have chosen $\uomega^{(k+1)}$. We  still have the ambiguity in $\ulambda^{(k+1)}$ up to constants along the fibres. We fix this by decreeing that the fibres have volume $1+O(\epsilon^{k+2})$ with respect to the
form $\umu^{[k+1]}$ defined by $\ulambda^{[k+1]}, \uomega^{[k+1]}$.

\

{\bf Step 4} Now we choose $H^{[k+1]}$ so that 
$  d_{H^{[k+1]}}\uomega^{[k+1]}$ and $  d_{H^{[k+1]}}\umu^{[k+1]}$  are $ O(\epsilon^{k+2})$. 

This completes the proof of the first item in Proposition 3.

For the second item, we take the power solution above for an input map $h$ and write

$$d_{H^{\epsilon}} \Theta^{\epsilon}= M(\epsilon, h)= \left(M_{0}(h) + \epsilon M_{1}(h) + \epsilon^{2} M_{2}(h)+\dots\right) \ulambda, $$
where  $M_{i}(h)$ are sections of the normal bundle of the image of $h$ in $H^{2}(X)$.  We know that $M_{0}(h)$ is the mean curvature $m(h)$,  whose linearisation with respect to normal variations is a Jacobi operator $J_{h}$ say. Locally, working over a ball, we can find a right inverse $G_{h}$ to the Jacobi operator. Then starting with a solution $h_{0}$ of the equation $m(h_{0})=0$ we can construct a formal power series solution 
$h^{\epsilon}= h_{0}+\sum \xi_{k} \epsilon^{k}$  to the equation $M(\epsilon, h^{\epsilon})=0$ in a standard fashion. For example the first term is given by
$$ h^{\epsilon}= h_{0} -\epsilon G_{h_{0}}(M_{1}(h_{0})) + O(\epsilon^{2}). $$
   
\subsection{The Torelli Theorem for K3 surfaces}

The question of which compact 4-manifolds admit hypersymplectic structures is an interesting open problem. But for hyperk\"ahler structures there is a complete classification. The underlying manifold must be a K3 surface or a 4-torus and we restrict here to the former case. We recall the fundamental {\it global Torelli Theorem}, for which see \cite{kn:BPV}, Chapter VIII for example. 
 \begin{prop}
     Let $X$ be the oriented 4-manifold underlying a K3 surface.  Write ${\cal C}$ for the set of classes $\delta\in H^{2}(X,\bZ)$ with $\delta.\delta=-2$ and write ${\rm Diff}_{0}$ for the group of diffeomorphisms of $X$ which act trivially on $H^{2}(X)$.  A hyperk\"ahler structure on $X$ determines  a self-dual subspace $H^{+}\subset H^{2}(X,\bR)$ (the span of $[\omega_{i}]$) which has the two properties:
\begin{itemize}
\item $H^{+}$ is a maximal positive subspace for the  cup product form.
\item  There is no class $\delta\in {\cal C}$ orthogonal to $H^{+}$.
\end{itemize}
Conversely, given a subspace $H^{+}\subset H^{2}(X;\bR)$ satisfying these two conditions there is a hyperk\"ahler structure with volume $1$ realising $H^{+}$ as the self-dual subspace and the structure is unique up to the action of ${\rm Diff}_{0}\times SL(3,\bR)$.
\end{prop}
(In the last sentence, the hyperk\"ahler structure is regarded as an element of $\bR^{3}\otimes \Omega^{2}(X)$ and $SL(3,\bR)$ acts on $\bR^{3}$.)

  \
  
 Let $h$ be a positive map from the $3$-ball to $H^{2}(X;\bR)$. We say that $h$ {\it avoids -2 classes} if for each point in $B$ there is no $\delta$ in ${\cal C}$ orthogonal to the image of the derivative. This is an open condition on $h$. The definition extends in the obvious way to the sections of the flat bundle ${\cal H}_{\chi}$ defined by a fibration $X\rightarrow M\rightarrow B$ as considered above. Then a hyperk\"ahler element $\uomega$
 defines a positive section of ${\cal H}_{\chi}$ which avoids $-2$ classes. A complete discussion of the converse  brings in topological questions about the group ${\rm Diff}_{0}$ which we do not want to go into here.  At least in the local situation, when $B$ is a ball, we can say that a positive section which avoids $-2$ classes defines a hyperk\"ahler element.
\section{Kovalev-Lefschetz fibrations}

Baraglia \cite{kn:B} shows that there are essentially no  co-associative fibrations over a compact base manifold $B$, so to get an interesting global problem we must  allow singularities. The singularities that we consider here are those  introduced by Kovalev \cite{kn:K}, where the fibres are allowed to develop ordinary double points.
\begin{defn} A  differentiable Kovalev Lefschetz (KL) fibration consists of data
$(M,B,L,\tilde{L},\pi)$ of the following form.
\begin{itemize} \item $M$ is a smooth oriented compact $7$-manifold, $B$ is a smooth oriented compact $3$-manifold and $\pi:M\rightarrow B$ is a smooth map.
\item $L\subset B$ is a  $1$-dimensional submanifold (i.e a link).
\item $\tilde{L}\subset M$ is a $1$-dimensional submanifold and the restriction of $\pi$ gives a diffeomorphism  from $\tilde{L}$ to $L$.
\item $\pi$ is a submersion outside $\tilde{L}$ and for any point $p\in \tilde{L}$ there are oriented charts around $p,\pi(p)$ in which $\pi$ is represented by a map
$$  (z_{1}, z_{2}, z_{3},  t)\mapsto (z_{1}^{2}+z_{2}^{2}+z_{3}^{2}+E(\underline{z},t), t), $$
from $\bR^{7}=\bC^{3}\times \bR$ to $\bR^{3}=\bC\times\bR$, where  $E(\underline{z},t)$ and its first and second partial derivatives vanish at $\underline{z}=0$.
\end{itemize}
\end{defn}

Thus the fibre $\pi^{-1}(p)$ of a KL fibration is a smooth $4$-manifold if $p$ is not in $L$, and if $p$ is in $L$ there is just one singular point, modelled on the quadric cone $\{z_{1}^{2}+z_{2}^{2}+z_{3}^{2}=0\}$ in $\bC^{3}$, which is homeomorphic to the quotient $\bC^{2}/\pm 1$. Now we define a torsion-free (respectively closed) co-associative KL fibration to be a  differentiable one, as above,  together with a smooth positive $3$-form $\phi$ on $M$ which is torsion-free (respectively closed) and such that the non-singular fibres are co-associative. We also require the orientations to be compatible with the form $\phi$.

A good reason for focusing on these kind of singularities is that Kovalev has outlined a programme that would produce  many examples of torsion-free co-associative fibrations of this shape \cite{kn:K}. (Kovalev's programme seems geometrically plausible, but there are some analytical difficulties in the available 2005 preprint \cite{kn:K}).   One can certainly consider other kinds of singularities, but we will not do so here.

 Note that a compact manifold with holonomy the full group $G_{2}$ has finite fundamental group \cite{kn:J}. Thus,  up to finite coverings (and using the solution of the Poincar\'e conjecture), the essential case of interest is when the base $B$ is the 3-sphere.

We can now adapt our previous discussion (for genuine fibrations) to the KL situation. We will only consider the case when the smooth fibres are diffeomorphic to a $K3$ surface.  Let $\pi:M\rightarrow B$ be a differentiable KL fibration and $\gamma$ be a small loop around one of the components of $L$. Then we can consider the {\it monodromy} around $\gamma$ of the co-homology of the fibre. It is a well-known fact that this monodromy is the reflection in a \lq\lq vanishing cycle' $\delta$ which is a class in $H^{2}(X;\bZ)$ with $\delta.\delta=-2$. (In particular, this monodromy  has order $2$, so we do not need to discuss the orientation of $\gamma$.)

To build in this monodromy we consider $B$ as an orbifold, modelled on
$(\bR\times \bC)/\sigma$ at points of $L$, where $\sigma$ acts as $-1$ on $\bC$ and $+1$ on $\bR$.  Recall that an orbifold atlas is a covering of $B$ by charts $U_{\alpha}\subset B$ such that
if $U_{\alpha}$ intersects $L$ it is provided with an identification with  a quotient $\tilde{U}_{\alpha}/\sigma$ for a $\sigma$-invariant open set $\tilde{U}_{\alpha}\subset \bR\times\bC$. A  {\it flat orbifold vector bundle over $B$} with fibre a vector space $V$ can be defined with reference to such an atlas as having local trivialisations of the form $U_{\alpha}\times V$ over the charts which do not meet $L$ and $(\tilde{U}_{\alpha}\times V, \tau_{\alpha})$ over those which do, where $\tau_{\alpha}$ is an involution of $V$. These trivialisations satisfy obvious compatibility conditions on the overlaps of charts. In  our situation the cohomology of the fibres provides such a flat orbifold vector bundle ${\cal H}$ over $B$ with the involutions $\tau_{\alpha}$ given by the reflections in the vanishing cycles. We also want to consider {\it flat orbifold affine bundles}. These are defined in just the same way except that  the $\tau_{\alpha}$ are replaced by affine involutions of $V$

Let ${\cal V}$ be any flat orbifold vector bundle over $B$. We can define a sheaf $\underline{{\cal V}}$ of locally constant sections of ${\cal V}$. In a chart $U_{\alpha}$ which intersect $L$ such a section is, by definition, given by a $\tau_{\alpha}$-invariant element of $V$. Any flat orbifold affine bundle determines a flat orbifold vector bundle via the natural map from the affine group of $\bR^{n}$ to $GL(n,\bR)$. We will say that the affine bundle is a {\it lift} of the vector bundle. Extending the standard theory for flat bundles recalled in Section 2 above; if 
${\cal V}$ is a flat orbifold vector bundle over $B$ then  there is a natural 1-1 correspondence between equivalence classes of flat orbifold affine bundles lifting ${\cal V}$ and elements of the cohomology group $H^{1}(B;\underline{{\cal V}})$. 

Write ${\cal V}_{\chi}$ for the affine bundle determined by   a flat orbifold vector bundle ${\cal V}$ and a class $\chi\in H^{1}(\underline{{\cal V}})$. We define the sheaf of smooth sections of ${\cal V}_{\chi}$ in  the standard fashion; given over  a chart $(\tilde{U}_{\alpha}, \tau_{\alpha})$  by  smooth maps 
$  f_{\alpha}: \tilde{U}_{\alpha}\rightarrow V $ satisfying the equivariance condition

$$f_{\alpha}(\sigma(x))=\tau_{\alpha}(f_{\alpha}(x)). $$

 We will want to consider {\it orthogonal} vector and affine bundles. That is, we suppose that the model vector space $V$ has a nondegenerate quadratic form of signature $(p,q)$ and all the bundle data, and in particular each involution $\tau_{\alpha}$, is compatible with this. For our purposes we restrict to the case when $p=3$. Away from $L$ our bundle ${\cal V}_{\chi}$ is an ordinary flat affine bundle so we are in the situation considered before and we have a notion of a {\it positive section}--i.e. the image of the derivative is a maximal positive subspace of the fibre. We want to extend this  notion over $L$. To discuss this we will now fix attention on the case when each $\tau_{\alpha}$ is a reflection in a negative vector. That is we assume that 
$$  \tau_{\alpha} (v)= v + (\delta_{\alpha} ,v) \delta_{\alpha}+ \lambda_{\alpha}\delta_{\alpha}, $$ where $ \delta_{\alpha}\in V$ with $(\delta_{\alpha},\delta_{\alpha})=-2$ and $\lambda_{\alpha}\in \bR$. Let $f_{\alpha}:\tilde{U}_{\alpha}\rightarrow V$ be an equivariant map, where $\tilde{U}_{\alpha}$ is some neighbourhood of the origin in $\bR\times \bC$. Without real loss of generality, and to simplify notation, we can suppose that $\lambda_{\alpha}=0$ and that $f_{\alpha}(0)=0$. We write $Df_{\alpha}, D^{2} f_{\alpha}$ for the first and second derivatives at $0$.
\begin{defn}
We say that the $f_{\alpha}$ is a {\em branched positive map} at $(0,0)$  if
\begin{enumerate} \item $Df_{\alpha}$ vanishes on $\bC$, so the image of $DF_{\alpha}$ is spanned by a single vector $v_{0}\in V$.
\item  The restriction of $D^{2}f_{\alpha}$ to $\bC$ can be written as $v_{1} {\rm Re} z^{2} + v_{2} {\rm Im} z^{2}$ where $v_{1}, v_{2}\in V$.
\item The vectors $v_{0}, v_{1}, v_{2}$ span  a maximal positive subspace in $\delta^{\perp}\subset V$. 
\end{enumerate}\end{defn}

In terms of local co-ordinates $(w,t)$  on $B$, where $w=z^{2}$, this can be expressed by saying that we have multi-valued function
$$  (w,t)\mapsto  v_{0} t+ v_{1} {\rm Re} w + v_{2} {\rm Im} w + O(w^{3/2}, t^{2}), $$
with an ambiguity of sign in in the component along $\delta\in V$. 

With this definition in place we can define a {\it  positive section} $h$ of our flat affine orbifold bundle ${\cal V}_{\chi}$ to be a section which is  positive away from $L$ and represented by a branched positive map  at points of $L$. We can define the 3-volume ${\rm Vol}_{3}(h)$  of the section in the obvious way. Finally we can define a {\it maximal  positive section} to be  a positive section which is maximal in the previous sense away from $L$. These definitions imply that near points of $L$ such a section yields a branched solution of the maximal submanifold equation, in the usual sense of the literature.

\subsection{Proposed adiabatic limit problem}

 Return now to a KL fibration $\pi:M\rightarrow B$,  the flat orbifold vector bundle ${\cal H}$ over $B$ and the corresponding sheaf $\underline{{\cal H}}$. 
It is straightforward to check that $\underline{\cal H}$ can be identified with the second derived direct image of the constant sheaf $\bR$ on $X$ and the first derived direct image is zero. So the Leray spectral sequence gives an exact sequence
     $$   0\rightarrow H^{3}(B;\bR)\rightarrow H^{3}(M;\bR)\stackrel{b}{\rightarrow} H^{1}(B;\underline{{\cal H}})\rightarrow 0. $$

\

\begin{prop}
Suppose that $\phi_{0}$ is a closed positive 3-form on $M$ with respect to which $\pi:M\rightarrow B$ is a co-associative KL fibration and  set $\chi= b([\phi_{0}])\in H^{1}((B;\underline{{\cal
H}})$. Then $\phi_{0}$ determines a branched positive section of ${\cal H}_{\chi}$.
\end{prop} 

The proof is not hard but we will not go into it here.

\

\

 We move on to the central point of this article. We say that a section $h$ of a flat affine bundle ${\cal H}_{\chi}$ {\it avoids excess $-2$ classes} if the following two conditions hold.
\begin{itemize}
\item The restriction to $B\setminus L$ avoids $-2$ classes, in the sense  defined before.
 \item At a point of $L$, and using the representation as in Definition 2 above,  the only elements of $ {\cal C}$ orthogonal to $v_{0}, v_{1}, v_{2}$ are $\pm \delta_{\alpha}$.
 \end{itemize}
 
  Write $[B]\in H^{3}(X)$ for the pull-back of the fundamental class of $B$. We makes a two-part conjecture.
\begin{conj}  Fix $[\phi_{0}]\in H^{3}(M)$ and $\chi=b([\phi_{0}])$.
\begin{enumerate} \item If there is a positive section $h$  of ${\cal H}_{\chi}$ which avoids excess $-2$ classes, then for sufficiently large $R$ there is a closed  positive 3-form on $M$ in the class $[\phi_{0}]+ R [B]$ with respect to which $\pi:M\rightarrow B$ has co-associative fibres. 
\item If, in addition, the section $h$ can be chosen to be maximal then, for sufficiently large $R$,  we can choose the positive 3-form to define a torsion free $G_{2}$-structure. \end{enumerate}
\end{conj}

These are precise conjectures which express part of a more general---but more vague---idea that for large $R$ the whole 7-dimensional theory should be expressed in terms of   positive  sections of ${\cal H}_{\chi}$. For example we would hope that there be should be converse statements, going from the existence of the structures on $M$, for large $R$,  to sections of ${\cal  H}_{\chi}$.  

We can extend the discussion to Bryant's Laplacian flow for closed positive 3-forms
$$  \frac{\partial \phi}{\partial t}= d d^{*}\phi. $$
The discussion in Section 2 suggests that the adiabatic limit should be the mean curvature flow for positive sections
\begin{equation}   \frac{\partial h}{\partial t}= m(h), \end{equation}
just as Hitchin's volume functional relates to the volume functional ${\rm Vol}_{3}{h}$.

Proving precise statements such as in Conjecture 1 will clearly involve substantial analysis. But we can leave that aside and formulate questions,  within the framework of sections of flat orbifold bundles, which can be studied independently.  Here we can take any orthogonal flat affine orbifold bundle ${\cal V}_{\chi}$ with orbifold structure defined by reflections, as above. 

\begin{enumerate}
\item 
 Does ${\cal V}_{\chi}$ admit positive sections?
 \item If so, does it admit a maximal positive section?
 \item If ${\cal V}_{\chi}$ admits positive sections, is the volume functional ${\rm Vol}_{3}$ bounded above?
\item Does the mean curvature flow (22) exist and converge to a maximal positive section?\end{enumerate}
Naturally we hope that these questions should have some relation with the corresponding questions for positive 3-forms at the beginning of this article. 
\

All of these questions should be viewed as tentative, provisional, formulations. We include two remarks in that direction.
\begin{itemize}
\item The condition of \lq\lq avoiding excess $-2$ classes'' seems very unnatural from the point of view of sections of the flat orbifold bundle ${\cal H}_{\chi}$. One could envisage that there could be cases where there are positive sections of ${\cal H}_{\chi}$ which avoid excess -2 classes but where there is a maximal positive section which does not satisfy this condition. This might correspond to a torsion-free $G_{2}$ structure on a singular space obtained by collapsing subsets of $M$. 
\item In the definition of a Kovalev-Lefschetz fibration we required that there be at most one singular point in each fibre and the critical values form a union of disjoint circles. This should be the \lq\lq generic'' situation, but one could envisage that in generic 1-parameter families one should allow two critical points in a fibre. In that vein one could imagine that the long-time definition of the mean curvature flow (22) would require modifying the link $L$ by allowing two components to cross, under suitable conditions on the monodromy of the flat bundle. Thus the right notion would not be a single KL fibration but an equivalence class generated by these operations.
\end{itemize}

\section{More adiabatic analogues}

\subsection{Negative definite cup product}
Let $M$ be any compact 7-manifold with a torsion-free $G_{2}$ structure defined by a form $\phi$. Suppose also that the holonomy is the full group $G_{2}$, which implies that $H^{1}(M;\bR)=0$ \cite{kn:J}. Then it is known via Hodge Theory that the quadratic form on $H^{2}(M;\bR)$ defined by
$$  Q_{[\phi]}( \alpha)=  \langle \alpha \cup\alpha\cup[\phi], [M]\rangle, $$
is negative definite \cite{kn:J}. Now suppose that $M$ has a $C^{\infty}$ Kovalev-Lefschetz fibration $\pi:M\rightarrow B$ with $K3$ fibres. The vanishing of $H^{1}(M)$ requires that $H^{1}(B;\bR)=0$. As in the previous section, let $[\phi_{0}]$ be any fixed class in $H^{3}(M)$ and $\phi_{R}= [\phi_{0}]+R [B]$. Then the limit of $R^{-1} Q_{[\phi_{R}]}$ as $R$ tends to infinity is given by restricting to a smooth fibre $X$ and taking the cup product on the fibre. Poincar\'e duality on $B$ gives $H^{2}(B;\bR)=0$ and then the Leray spectral sequence shows that $H^{2}(M;\bR)= H^{0}( \underline{{\cal H}})$, which can be identified with the subspace of $H^{2}(X)$  preserved by the monodromy of the fibration. 

Now we have a related result in the context of positive sections.
\begin{prop}
Let ${\cal V}$ be a flat orbifold bundle over a compact base $B$ and ${\cal V}_{\chi}$ be a lift to an affine orbifold bundle. Suppose that ${\cal V}_{\chi}$ admits a positive section and that  $H^{1}(B, \bR)=0$. If $\sigma$ is a non-trivial global section of $\underline{{\cal V}}$ then   $\sigma.\sigma<0$ at each point of $B$.
\end{prop}
We can write ${\cal V}= \bR\ \sigma \oplus {\cal V'}$. The condition that $H^{1}(B;\bR)=0$ implies that $\chi$ is trivial in the $\bR\ \sigma$ component, so  if $h$ is the positive section of ${\cal V}_{\chi}$ then there is a well-defined continuous function $\sigma.s$ on $B$ which attains a maximum at some point $p\in B$. Consider first the case when $p$ is not in $L$. The maximality implies  that $\sigma(p)$ is orthogonal to the image of $dh$ at $p$. By the definition of a positive section, this image is a maximal positive subspace which means that $\sigma(p).\sigma(p)<0$. The definition of a positive section at points of $L$ means that the same argument works when $p\in L$.

\subsection{Associative submanifolds}

Recall that a 3-dimensional submanifold $P\subset M$ in a 7-manifold with a positive $3$-form $\phi$ is called {\it associative} if the restriction of $\phi$ to the orthogonal complement of $TP$ in $TM$ vanishes at each point of $P$. This implies that, with the right choice of orientation, the restriction of $\phi$ to $TP$ is the volume form. These associative submanifolds are further examples of Harvey and Lawson's calibrated submanifolds \cite{kn:HL}, with calibrating form $phi$ and there is an elliptic deformation theory \cite{kn:McL}.   They are interesting geometrical objects to investigate in  $G_{2}$-manifolds:  in particular there is a potential connection with moduli problems, because if a non-zero homology class $\pi \in H_{3}(M)$ can be represented by a compact associative submanifold then we have $\langle [\phi], \pi \rangle>0$. 

Let $X$ be a  K3 surface with a hyperk\"ahler triple $\omega_{1}, \omega_{2}, \omega_{3}$. The cohomology classes $[\omega_{i}]$ span a positive subspace $H^{+}\subset H^{2}(X;\bR)$. Each non-zero vector $\tau=(\tau_{1}, \tau_{2}, \tau_{3})$ defines a complex structure $I_{\tau}$ on $X$ and a symplectic form $\omega_{\tau}=\sum \tau_{i}\omega_{i}$. Write $H^{\perp}_{\tau}$ for the orthogonal complement of $[\omega_{\tau}]$ in $H^{+}$.  If $\Sigma$ is a smooth complex curve in $X$, for the complex structure $I_{\tau}$, then the product of $\Sigma$ and the line $R\tau$ is an associative submanifold in $X\times \bR^{3}$ with the  positive $3$-form $\sum \omega_{i} dt_{i} - dt_{1}dt_{3} dt_{3}$. In particular, suppose that $c\in H^{2}(X;\bZ)$ is a class with $c^{2}=-2$ and $c$ is orthogonal $H^{\perp}_{\tau}$; then standard theory tells us that there is a unique such curve $\Sigma$ in the class $c$, and $\Sigma$ is a $2$-sphere. This leads to a natural candidate for an \lq\lq adiabatic description'' of certain associative submanifolds in a Kovalev-Lefschetz fibration. We begin with the local picture, away from the critical set $L\subset B$. Suppose that we have a positive section given by a map $h$ into $H^{2}(X;\bR)$ and fix a class $c$ as above. Let $\gamma(s)$ be an embedded path in $B$ so $h\circ\gamma$ is a path in $H^{2}(X;\bR)$ and for each $s$ the derivative of $h\circ\gamma$ defines a complex structure on $X$ and a $2$-dimensional subspace $H^{\perp}(s)\subset {\rm Im} \ dh \subset H^{2}(X;\bR)$.  If $c$ is orthogonal to $H^{\perp}(s)$ for each $s$ then we have a complex $2$-sphere in the fibre over $\gamma(s)$ and the union of these 2-spheres, as $s$ varies, is our candidate for an adiabatic limit of associative submanifolds. The condition on the derivative of $h\circ\gamma$ has a simple geometric meaning. Up to parametrisation of $\gamma$, it is just the condition that $\gamma$ be a {\it gradient path}  for the function $\underline{t}\mapsto c. h(\underline{t})$ with respect to the metric on $B$ induced by $h$. 

To bring in the critical set $L\subset B$, recall that if $p\in B\setminus L$ is a point close to $L$ there is a vanishing cycle $\Sigma$ in the fibre $X_{p}$ over $p$ which is a $2$-sphere bounding a \lq\lq Lefschetz thimble'' in $M$, and the homology class of the vanishing cycle is a class $c$ as above. Let $\gamma:[0,1]\rightarrow B$ be a path with $\gamma(0)$ and $ \gamma(1)$ in  $L$ and $\gamma(s)\in B\setminus L$ for $s\in (0,1)$. For small $\epsilon>0$ there are vanishing cycles (defined up to isotopy) in the fibres over $\gamma(\epsilon), \gamma(1-\epsilon)$ and we say that $\gamma$ is a {\it $C^{\infty}$ matching path} if these agree under parallel transport  along $\gamma$. In this case we can construct an embedded 3-sphere in $M$ by capping off the 2-sphere fibration over $(\epsilon, 1-\epsilon)$ with the two Lefschetz thimbles.  We say that $\gamma$ is a {\it homological matching path} if the homology classes of the two vanishing cycles   agree, a condition which is determined by the flat orbifold bundle ${\cal H}$. For the present discussion we will not  distinguish between the homological and $C^{\infty}$ conditions.
Now we say that $\gamma$ is a {\it matching gradient path} if it is a homological matching path and if over the interior $(0,1)$ it is a gradient path of the function $h.c$ (interpreted using parallel transport of $c$ along $\gamma$). These are plausible candidates for the adiabatic limits of certain associative spheres. More precisely, one can hope that, provided $\gamma$ is a \lq\lq nondegenerate'' solution of the matching gradient path condition in a suitable sense, then there will be  an associative 3-sphere in the $G_{2}$ manifold for large values of $R$. On the other hand, the matching gradient paths are described entirely in terms of ${\cal H}_{\chi}$. The definition extends immediately to any flat orbifold bundle and they  are objects which can be studied independent of any connection with $G_{2}$ geometry. This proposed description of associative submanifolds  is in the same vein as constructions of Lagrangian spheres in symplectic topology \cite{kn:S} and of tropical descriptions of holomorphic curves in Calabi-Yau manifolds \cite{kn:G}.

\section{Variants}\subsection{Torus fibres, special Lagrangian fibrations}
The local discussion of co-associative $G_{2}$-fibrations in Section 2 applies with little change when the fibres are 4-tori rather than $K3$ surface. The \lq \lq adiabatic limit data'' is for a maximal submanifold in $\bR^{3,3}=H^{2}(T^{4})$. The difference in this case is that such a maximal submanifold yields a genuine exact solutions of the equations for non-zero $\epsilon$.  This was the case treated by Baraglia in \cite{kn:B}. The source of this difference is that the curvature of the connection $H$ vanishes in this situation---the horizontal subspaces are tangent to sections of the fibration. 

Suppose now that $M=S^{1}\times N$ for a $6$-manifold $N$, and consider $G_{2}$-structures on $M$ induced from $SU(3)$ structures on $N$. Then we have
$$  \phi= \rho +\Omega d\theta, $$
where $d\theta$ is the standard $1$-form on $S^{1}$ and $\rho,\Omega$ are respectively a $3$-form and a $2$-form on $N$.  Here we consider the case of  torus fibrations $M\rightarrow B$  induced from {\it special Lagrangian fibrations} $\pi_{N}:N\rightarrow B$ by taking a product of the fibres with a circle.  Now $$  H^{2}(T^{3}\times S^{1}) = H^{2}(T^{3})\oplus H^{1}(T^{3})$$ and the quadratic form on $H^{2}(T^{3}\times S^{1})$ is induced from  the dual pairing between $H^{1}(T^{3}), H^{2}(T^{3})$. In other words the reduction to the product situation yields a natural splitting
$$  \bR^{3,3}= \bR^{3}\oplus \left(\bR^{3}\right)^{*}. $$
The dual pairing also defines  a standard symplectic form on $\bR^{3,3}$ and the product structure implies that the image of the map $h:B\rightarrow \bR^{3,3}$ is Lagrangian. It follows that this the image is   the graph of the derivative of a function $F$ on   $\bR^{3}$. The basic fact, observed by Hitchin \cite{kn:H2}, is that a Lagrangian submanifold of $\bR^{3,3}$  defined as the graph of the derivative of a function $F$  is a maximal submanifold if and only if $F$ satisfies the Monge-Amp\`ere equation 
$$  {\rm det} \left( \frac{\partial^{2}F}{\partial x_{i}\partial x_{j}}\right) = C$$ for a  constant $C$.
This makes a bridge between the co-associative discussion and the extensive literature on special Lagrangian fibrations and \lq\lq large complex structure'' limits for Calabi-Yau 3-folds. 
\subsection{K3 fibrations of Calabi-Yau $3$-folds}

Suppose again that $M=S^{1} \times N$ but now that $B=S^{1}\times S^{2}$ and consider fibrations which are products with $S^{1}$ of $\varpi_{N}: N\rightarrow S^{2}$ with $K3$ fibres. More precisely, we suppose that $\varpi_{N}$ is a Lefschetz fibration in the usual sense, with a finite set of critical values in $S^{2}$. The class $[\omega]$ gives a fixed class in the cohomology  of the fibres so now we form a flat orbifold vector bundle  ${\cal H}$ with fibre $R^{2,19}$. The discussion of positive sections goes over immediately to this 2-dimensional case. We could hope that this has a bearing on the question: when can a differentiable Lefchetz fibration with $K3$ fibres be given the structure of a holomorphic fibration of a Calabi-Yau 3-fold?

We can again make a bridge between this point of view and other, more standard, ones. Recall that the Grassmann manifold ${\rm Gr}^{+}$ of positive oriented 2-planes in $\bR^{2,19}$ can be identified with an open subset of a quadric $Q\subset \bC\bP^{20}$. This is achieved by mapping the plane spanned by the orthonormal frame $e_{1}, e_{2}$ to the complex line $[e_{1} + i e_{2}]\in \bC\bP^{20}$. 
Let $\Sigma\subset \bR^{2,19}$ be a surface with positive tangent space at each point. Then we have a Gauss map from $\Sigma$ to ${\rm Gr}^{+}$. Just as in the Euclidean case, the maximal condition implies that this is holomorphic with respect to the induced complex structure on $\Sigma$; hence the  image is a complex curve. Conversely, start with a complex curve in ${\rm Gr}^{+}\subset \bC\bP^{20}$, lift it locally to $\bC^{21}$ and parametrise it by a complex variable $z\in U\subset \bC$. So we have holomorphic functions $f_{0}(z), \dots, f_{20}(z)$ with $Q(f_{0}, \dots ,f_{20})=0$. Let $H_{i}$ be the integrals, solving
$\frac{dH_{i}}{dz}= f_{i}$. Then the map from $U$ to $\bR^{2,19}$ with components ${\rm Re} H_{i}$ has image a maximal surface in $\bR^{2,19}$. This is the usual Weierstrasse construction. In our context, the complex curve in ${\rm Gr}^{+}$ amounts to a holomorphic description of a family of polarised $K3$ surfaces in terms of the Hodge data $H^{2,0}\subset H^{2}(X,\bC)$ and the Weierstrasse  construction makes the bridge with the maximal submanifold equation.

\subsection{Cayley fibrations}

Here we consider the local adiabatic problem for a fibration of an 8-dimensional manifold with holonomy ${\rm Spin}(7)$ by Cayley submanifolds. Recall that a torsion-free ${\rm Spin}(7)$ structure can be defined by a closed $4$-form which is algebraically special at each point \cite{kn:J}.  As the basic model we take two oriented $4$-dimensional Euclidean spaces $V_{1}, V_{2}$ and a fixed isomorphism of the respective  self-dual forms $\Lambda^{2}_{+}(V_{i}^{*})$. Then we have a 4-form
$$  \Omega_{0}= d{\rm vol}_{1} + \sum \omega_{a}\wedge \omega_{a}'+ d{\rm vol}_{2} $$ on $V_{1}\oplus V_{2}$, where $d{\rm vol}_{i}$ are the volume forms in the two factors, $\omega_{a}$ runs over a standard orthonormal basis of $\Lambda^{2}_{+}(V_{1}^{*})$ and $\omega_{a}'$ the corresponding basis for $\Lambda^{2}_{+}(V_{2}^{*})$. A closed $4$-form $\Omega$ on an $8$-manifold
$M$ which is algebraically equivalent to $\Omega_{0}$ at each point defines a torsion free ${\rm Spin}(7)$ structure. A 4-dimensional submanifold $X\subset M$ is a Cayley submanifold if at each point of $X$ the triple $(TM_{x}, \Omega(x), TX_{x})$ is algebraically equivalent to the model triple $(V_{1}\oplus V_{2}, \Omega_{0}, V_{1})$. These are, again, examples of calibrated submanifolds \cite{kn:HL} and there is an elliptic deformation theory \cite{kn:McL}.

A fibration  $\pi:M\rightarrow B$ with Cayley fibres can be expressed in the manner of Section 2. There is a  connection $H$ which furnishes a decomposition of the forms on $M$ into $(p,q)$ components, and a closed  4-form $\Omega$ with three components
$$ \Omega = \Omega_{4, 0} + \Omega_{2,2}+ \Omega_{0,4}$$
The algebraic data is that at each point $x$ of $M$ we have oriented Euclidean structures on the horizontal and vertical subspaces, and an isomorphism between the spaces of self-dual forms, as above. The condition that $\Omega$ is closed can be viewed a set of coupled equations for the connection and the three components. A discussion parallel to that in the $7$-dimensional case suggests that in the adiabatic limit the equations should de-couple,  forcing the structure on the each fibre to be hyperk\"ahler.  This analysis leads us to formulate the following definition.
\begin{defn} Let ${\cal H}$ be a real vector space with a quadratic form of signature $(3,q)$. {\em Adiabatic Cayley data} for ${\cal H}$ consists of 
an oriented Riemannian $4$-manifold $(B,g)$ and a closed ${\cal H}$-valued self-dual $2$-form $\Psi$ on $B$ such that $\Psi^{*}\Psi$ is the identity endomorphism of  $\Lambda^{2}_{+}$.
\end{defn}

To explain the notation here, at each point $b\in B$, the value $\Psi(b)$ is regarded as a linear map from the fibre of $\Lambda^{2}_{+}$ to ${\cal H}$ (using the metric on $\Lambda^{2}_{+}$). Then $\Psi^{*}(b)$ is the adjoint, defined using the quadratic forms on each space, so the composite $\Psi^{*}\Psi$ is an endomorphism of $\Lambda^{2}_{+}$. For our immediate application we take ${\cal H}$ to be the 2-dimensional cohomology of a $4$-torus or a $K3$ surface $X$. Given adiabatic Cayley data as above, the image of  $\Psi(b)$ at each point $b$ is a maximal positive subspace of ${\cal H}$ so defines a hyperk\"ahler structure on $X$. (More precisely this is defined up to diffeomorphism, but that is not important for the local discussion here.) If $\omega_{i}, \omega_{2}, \omega_{3}$ is a standard orthonormal frame for $\Lambda^{2}_{+}(b)$ then we have a corresponding hyperk\"ahler triple  $\omega'_{i}$ of $2$-forms on the fibre $\pi^{-1}(b)$. That is, $$[\omega'_{i}]= \Psi(b)(\omega_{i})\in H^{2}(\pi^{-1}(b)). $$
The sum $$\Omega_{2,2}= \sum \omega_{a}\wedge\omega'_{a}$$ is independent of the choice of orthornormal frame and
$$ d_{f}\Omega_{2,2}=0. $$ Then is then a unique way to choose a connection $H$ which is compatible with the volume form on the fibres and so that
$$   d_{H} \Omega_{2,2}=0.$$
We take $\Omega_{0,4}$ to be extension of the volume form on the fibres by this connection and $\Omega_{4,0}$ to be the lift of the volume form on the base. In this way we define a $4$-form $\Omega$ on the total space of he correct algebraic type and with
    $$  d\Omega= F_{H} \Omega_{0,4}+ F_{H} \Omega_{2,2}, $$
    terms which, as we argued before, should be suppressed in the adiabiatic limit. (When the fibre is a 4-torus then, just as in Baraglia's case discussed above, the curvature $F_{H}$ is zero, so we can  write down genuine torsion-free ${\rm Spin}(7)$-structures, given adiabatic Cayley data. )

    One remark about Definition 3 is that the metric $g$ is entirely determined by the ${\cal H}$-valued form $\Psi$. In other words suppose we are given an oriented $4$-manifold $B$ and an ${\cal H}$ valued $2$-form $\Psi$ on $B$. Then we can define $\Psi^{*}\Psi$, up to ambiguity by a positive scalar, using the wedge product form on $2$-forms. Let us say that $\Psi$ is {\it  special} if at each point $\Psi^{*}\Psi$ is a projection onto a $3$-dimensional positive subspace in $\Lambda^{2}$ (up to a factor). This subspace defines a conformal structure on $B$---the unique conformal structure so that this subspace is self-dual---and we fix a metric by by normalising so that $\Psi^{*}\Psi$ is the identity on this subspace. In sum, an equivalent description of adiabatic Cayley data is a closed, special, ${\cal H}$-valued $2$-form.  
    
    \subsection{Ricci curvature}

    We recall first a standard result.
    \begin{prop}
    Let $V\subset \bR^{p,q}$ be a $p$-dimensional maximal positive submanifold. Then the Ricci curvature of the induced Riemannian metric on $V$ is non-negative.
\end{prop}
We review the proof. By elementary submanifold theory, the Riemann curvature tensor of $V$ is expressed in  terms of the second fundamental form $S$. Fix a point in $V$ and orthonormal frames for the tangent and normal bundles so $S$ has components $S_{ij}^{a}$ where $i,j$ label the tangent frame and $a$ the normal frame. Then the Riemann curvature tensor is
$$  R_{ij kl} = \sum_{a} S^{a}_{ik} S^{a}_{jl}- S^{a}_{jk} S^{a}_{il}, $$
where we have an opposite sign to the familiar Euclidean case because the quadratic form is negative in the normal direction. This gives
$$   R_{ik}= - \sum_{a,j} (S^{a}_{jj}) S^{a}_{ik}+ \sum_{a,j} S^{a}_{ij}S^{a}_{j k}. $$
The term $\sum_{j} S^{a}_{jj}$ is the mean curvature, which vanishes by assumption, and $S$ is symmetric so we have
$$   R_{ik}= \sum_{a,j} S^{a}_{ij}S^{a}_{kj}, $$ which is non-negative. More precisely, the Ricci curvature in a tangent direction $\xi$ to $V$ is the square of the norm of the derivative of the Gauss map along $\xi$.

\

The relevance of this for us is an immediate consequence.

\begin{cor} Suppose that ${\cal H}_{\chi}$ is a flat affine orbifold bundle over $(B,L)$ which admits a maximal positive section. Then the induced Riemannian metric on $B\setminus L$ has non-negative Ricci curvature. \end{cor}

We pass on to the higher dimensional case, where the result is less standard.
\begin{prop}
Suppose that $(B,g,\Psi)$ is a set of adiabatic Cayley data. Then the Ricci curvature of $g$ is non-negative. 
\end{prop}
We can regard $\Psi$ as giving an isometric embedding of the bundle $\Lambda^{2}_{+}$ over $B$ into the trivial bundle with fibre ${\cal H}$. There is thus a \lq \lq second fundamental form'' $S$ of this subbundle, or equivalently the derivative of the map from $B$ to the Grassmannian. Let ${\cal H}^{\perp}$ be the complementary bundle in ${\cal H}$. Then $S$ is a tensor in $\Lambda^{2}_{+}\otimes \Lambda^{1} \otimes {\cal H}^{\perp}$. A little thought shows that the condition that $d\Psi=0$ is equivalent to
\begin{itemize}
\item The connection on $\Lambda^{2}_{+}$ induced from the Levi-Civita connection is the same as the connection induced from the embedding $\Lambda^{2}_{+}\subset {\cal H}$;
\item The image of $S$ in $\Lambda^{3}\otimes {\cal H}^{\perp}$ induced the wedge product $\Lambda^{1}\otimes \Lambda^{2}_{+}\rightarrow \Lambda^{3}$ is zero.
\end{itemize}
  
The first item means that the Riemannian curvature of $\Lambda^{2}_{+}$ can be computed from $S$. Indeed it is given by $q_{{\cal H}^{\perp}}(S)$  where $q_{{\cal H}^{\perp}}$ is the quadratic map from $\Lambda^{+}\otimes \Lambda^{1}\otimes {\cal H}^{\perp}$ to $\Lambda^{2}\otimes \Lambda^{+}$ by the tensor product of
\begin{itemize}
\item The skew-symmetric wedge product $\Lambda^{1}\otimes \Lambda^{1}\rightarrow \Lambda^{2}$;
\item The skew-symmetric cross-product $\Lambda^{2}_{+}\otimes \Lambda^{2}_{+}\rightarrow \Lambda^{2}_{+}$;
\item The symmetric inner product ${\cal H}^{\perp}\otimes {\cal H}^{\perp}\rightarrow \bR$.\end{itemize}
Of course to pin down signs it is  crucial  that the last form is negative definite, but otherwise it will not  enter the picture. We define a quadratic map $q$ from $\Lambda^{+}\otimes \Lambda^{1}$ to $\Lambda^{2}\otimes \Lambda^{+}$ by using the first two components above. 

\

Next we want to recall how to compute the Ricci curvature of an oriented Riemannian 4-manifold $(B,g)$ from the curvature tensor $F$ of $\Lambda^{2}_{+}$. This curvature tensor $F$ lies in $(\Lambda_{+}\oplus\Lambda_{-})\otimes \Lambda_{+}$ so has two components, say
$$  F_{+}\in \Lambda_{+}\otimes\Lambda_{+}    \ \ , \ \ F_{-}\in \Lambda_{-}\otimes\Lambda_{+}. $$ We work at a fixed point $b\in B$ and take a unit tangent vector $e_{0}$. This choice of $e_{0}$ induces an isomorphism between the self-dual and anti-self-dual forms at $b$. To fix signs,  say that this isomorphism takes the self-dual form
$dx_{0}dx_{1}+ dx_{2} dx_{3}$ to the anti-self-dual form $-dx_{0}dx_{1}+ dx_{2} dx_{3}$. Using this isomorphism we can define a trace
$$  {\rm Tr}_{e_{0}}: \Lambda{-}\otimes \Lambda_{+}\rightarrow \bR, $$ depending on $e_{0}$. Of course we also have an invariant trace:
$$  {\rm Tr}: \Lambda_{+}\otimes \Lambda_{+}\rightarrow \bR. $$
We define $T_{e_{0}}:\Lambda^{2}\otimes \Lambda^{+}\rightarrow \bR $ to be the sum ${\rm Tr} + {\rm Tr}_{e_{0}}$ on the two factors. The formula we need is that
$$  {\rm Ric}(e_{0})= T_{e_{0}}(F). $$
This is easy to check, using the well-known decomposition of the curvature tensor of a 4-manifold, or otherwise. So in our situation
$$  {\rm Ric}(e_{0})= T_{e_{0}}(q_{{\cal H}^{\perp}}(S)). $$
So we need to show that for any $s\in \Lambda^{+}\otimes \Lambda^{1}$ which is in the kernel of the wedge product $\Lambda^{+}\otimes \Lambda^{1}\rightarrow \Lambda^{3}$ we have $ T_{e_{0}}(q(s))\leq 0$. Calculate in an orthonormal frame $e_{i}$ extending $e_{0}$. This induces a standard orthonormal frame $\omega_{i}$ for $\Lambda_{+}$. We can write
$s= \sum_{i=1}^{3} v_{i}\otimes \omega_{i}$,
where $v_{i}= t_{i} e_{0} + \sum_{i=1}^{3} s_{ij} e_{j}$. 
One finds that the condition that $s$ is in the kernel of the wedge product is that the the trace $\sum s_{ii}$ is zero and that the $t_{i}$ are determined by the skew symmetric part of  $s_{ij}$:
$$   t_{i} = s_{jk}- s_{kj}$$ for $i,j,k$ cyclic. On the other hand one finds that for any $s$ in $\Lambda^{1}\otimes \Lambda^{+}$ 
$$  T_{e_{0}}(q(s))= \sum_{{\rm cyclic}} t_{i}(s_{jk}-s_{kj}). $$
So for $s$ in the kernel of the wedge product this becomes $-\sum t_{i}^{2}$ which has the desired sign

In fact this calculation shows that again ${\rm Ric}(e_{0})$ is the square of the norm of the derivative of the map to the Grassmannian of $3$-planes in ${\cal H}$,   in the direction $e_{0}$.

\end{document}